\documentclass[a4paper, 11pt]{amsart}

\usepackage{amsthm,amsfonts}
\usepackage{amssymb,graphicx,color}
\usepackage{xurl}
\usepackage{hyperref}
\hypersetup{colorlinks=true,linkcolor=blue,citecolor=blue,urlcolor=blue}
\usepackage{enumitem}
\usepackage[margin=1.3in]{geometry}

\newtheorem{innercustomthm}{Theorem}
\newenvironment{customthm}[1]
  {\renewcommand\theinnercustomthm{#1}\innercustomthm}
  {\endinnercustomthm}
\newtheorem{innercustomcor}{Corollary}
\newenvironment{customcor}[1]
  {\renewcommand\theinnercustomcor{#1}\innercustomcor}
  {\endinnercustomcor}

\newtheorem{theorem}{Theorem}[section]
\newtheorem{lemma}[theorem]{Lemma}
\newtheorem{corollary}[theorem]{Corollary}
\newtheorem{proposition}[theorem]{Proposition}
\theoremstyle{definition}
\newtheorem{definition}[theorem]{Definition}
\newtheorem{example}[theorem]{Example}
\newtheorem{remark}[theorem]{Remark}
\theoremstyle{plain}
\newtheorem{question}{Question}
\newtheorem{conjecture}{Conjecture}

\newcommand{\R}{\mathbb{R}}
\newcommand{\C}{\mathbb{C}}
\newcommand{\Z}{\mathbb{Z}}
\newcommand{\Q}{\mathbb{Q}}
\newcommand{\Sp}{\mathbb{S}}
\newcommand{\Sing}[1]{{\rm Sing}(#1)}
\newcommand{\Reg}[1]{{\rm Reg}(#1)}

\begin{document}

\title[The metric Zariski multiplicity conjecture in dimension 4]
{The metric Zariski multiplicity conjecture in dimension 4}

\author[J. E. Sampaio]{Jos\'e Edson Sampaio}
\address{Jos\'e Edson Sampaio: Departamento de Matem\'atica, Universidade Federal do Cear\'a,
	      Rua Campus do Pici, s/n, Bloco 914, Pici, 60440-900,
	      Fortaleza-CE, Brazil.}
 \email{edsonsampaio@mat.ufc.br}



\thanks{The author was partially supported by CNPq-Brazil grants 303375/2025-6 and 404611/2025-7. This work was supported by the Serrapilheira Institute (grant number Serra -- R-2110-39576).}

\begin{abstract}
Recently, OpenAI announced a counterexample to the Zariski multiplicity conjecture in dimension 4. Here, we prove that the metric Zariski multiplicity conjecture is true in dimension 4. More precisely, let $f,g\colon(\mathbb{C}^n,0)\to(\mathbb{C},0)$ be germs of reduced holomorphic functions, and let $\varphi\colon(\mathbb{C}^n,0)\to(\mathbb{C}^n,0)$ be a bi-Lipschitz homeomorphism such that $\varphi(V(f))=V(g)$. We prove that if every irreducible component of the tangent cone of $V(f)$ at the origin has a singular locus of dimension at most two, then $V(f)$ and $V(g)$ have the same multiplicity at the origin. In particular, the multiplicity of complex analytic hypersurface singularities in $\mathbb{C}^n$, $n\leq 4$, is an embedded bi-Lipschitz invariant; the cases $n=2$ and $n=3$ were known. This extends a previous result of the author for tangent cones whose components have singular loci of dimension at most one. We also prove that two irreducible homogeneous hypersurfaces whose singular loci have dimension at most two, and which are mapped onto each other by a homeomorphism of $\mathbb{C}^n$ fixing the origin, have the same degree.
\end{abstract}

\maketitle
\section{Introduction}

Consider a reduced holomorphic function germ $f\colon(\C^n,0)\to(\C,0)$ with zero set $V(f)$. We can expand $f$ as $f=f_m+f_{m+1}+\cdots$, where each $f_k$ is a homogeneous polynomial of degree $k$ and $f_m\neq0$. Here, $m(V(f),0)=m$ denotes the \textbf{multiplicity} of $V(f)$ at $0$, and $in(f)=f_m$ is the \textbf{initial form} of $f$. In 1971, Zariski posed the following question \cite{Zariski:1971}:

\begin{question}[Zariski]\label{question:zariski}
Let $f,g\colon(\C^n,0)\to(\C,0)$ be germs of reduced holomorphic functions. If there is a homeomorphism $\varphi\colon(\C^n,0)\to(\C^n,0)$ such that $\varphi(V(f))=V(g)$, is it true that $m(V(f),0)=m(V(g),0)$?
\end{question}

This question, in its stated version, is famously known as the Zariski multiplicity conjecture. It has been proven true for $n=2$ \cite{Zariski:1932}, while the $n=3$ case remains unresolved (refer to \cite{Eyral:2007} for a comprehensive survey). Recently, OpenAI reported counterexamples to Question \ref{question:zariski}: in \cite{OAI26}, they constructed two holomorphic function germs on $(\C^4,0)$ with isolated singularities and multiplicities $4$ and $5$ that have ambiently homeomorphic zero sets. A related construction in \cite{OAI26b} yields germs with isolated singularities and multiplicities $2$ and $3$ exhibiting the same property in $\C^N$ (where $N$ is very large and a multiple of $8$). Should these results be verified, the Zariski multiplicity conjecture would be false for all dimensions $n\geq 4$. 

The metric counterpart of this conjecture is still open, and in light of OpenAI's examples, it might represent the more viable problem to tackle.

\begin{conjecture}[Metric Zariski multiplicity conjecture]\label{conj:metric}
Let $f,g\colon(\C^n,0)\to(\C,0)$ be germs of reduced holomorphic functions. If there is a bi-Lipschitz homeomorphism $\varphi\colon(\C^n,0)\to(\C^n,0)$ such that $\varphi(V(f))=V(g)$, then $m(V(f),0)=m(V(g),0)$.
\end{conjecture}

Conjecture \ref{conj:metric} is known to hold for $n\leq 3$ \cite[Theorem 2.3]{FernandesS:2016} (with the $n=2$ case already established by \cite{Zariski:1932}). It also holds when the irreducible components of the tangent cone of $V(f)$ have isolated singularities \cite[Theorem 2.2]{FernandesS:2016}, when their singular loci have dimension $\leq 1$ \cite[Corollary 5.7 and Remark 5.13]{Sampaio:2020} (see also \cite[Theorem 6.6]{FernandesS:2023}), and when $m(V(f),0)=2$ \cite{FernandesJS:2025}. If the bi-Lipschitz homeomorphism is not required to be ambient, multiplicity remains an invariant for complex analytic set germs of dimension $1$ (see \cite{P-T, F, N-P}) and dimension $2$ \cite{FernandesFS:2018}. However, this invariance fails for dimensions $\geq 3$ \cite{BirbrairFSV:2020, FernandesJS:2023} (consult \cite{FernandesS:2023} for an overview). In this work, we establish the following:

\begin{customthm}{\ref*{thm:main}}
Let $f,g\colon(\C^n,0)\to(\C,0)$ be germs of reduced holomorphic functions, and let $\varphi\colon(\C^n,0)\to(\C^n,0)$ be a bi-Lipschitz homeomorphism such that $\varphi(V(f))=V(g)$. If every irreducible component of the tangent cone $C(V(f),0)$ has a singular locus of dimension $\leq2$, then $m(V(f),0)=m(V(g),0)$.
\end{customthm}

Because any irreducible component of the tangent cone of a hypersurface in $\C^4$ is a $3$-dimensional cone whose singular locus has dimension at most $2$, this directly yields:

\begin{customcor}{\ref*{cor:n4}}
Conjecture \ref{conj:metric} has a positive answer for $n\leq4$.
\end{customcor}

Corollary \ref{cor:n4} stands in sharp contrast to the topological setting: if the findings in \cite{OAI26} hold, this corollary implies that no bi-Lipschitz homeomorphism of $(\C^4,0)$ can map their respective zero sets to one another (likewise, by \cite{FernandesJS:2025}, the homeomorphisms in \cite{OAI26b} cannot be bi-Lipschitz).

The novelty of Theorem \ref{thm:main} lies in handling tangent cones where components $V(P)$ possess $2$-dimensional singular loci---that is, the projective hypersurface $\mathbb{P} V(P)$ is singular along curves (e.g., $P=x^2z-y^2w$ on $\C^4$, as discussed in Example \ref{exam:umbrella}).

Our proof begins with a standard reduction: relying on \cite{Sampaio:2016} and \cite{FernandesS:2016}, Theorem \ref{thm:main} reduces to the case of homogeneous hypersurfaces. Specifically, it suffices to prove that any two irreducible homogeneous hypersurfaces $X=V(P)$ and $Y=V(Q)$ in $\C^n$ mapped to each other by a bi-Lipschitz homeomorphism $\Psi\colon\C^n\to\C^n$ have the same degree. Previous work \cite{FernandesS:2016, Sampaio:2020, FernandesS:2023} established this when the Milnor fibre $P^{-1}(1)$ has a non-vanishing Euler characteristic (see Proposition \ref{prop:lefschetz}), using L\^e's topological invariance of the monodromy \cite{Le:1973a} and the vanishing of the Lefschetz numbers of the monodromy iterates $h^k$ for $0<k<{\rm deg}(P)$. For the general case, we apply the Parusi\'nski--Pragacz formula \cite{ParusinskiP:1995} for the Euler characteristic of a singular projective hypersurface. When the singular locus of $\mathbb{P} V(P)$ has dimension $\leq 1$, every term in this formula---except the degree of $P$ itself---is invariant under ambient homeomorphisms mapping $V(P)$ to $V(Q)$. We control the degrees of the singular curves by analyzing the cohomology torsion of the corresponding punctured cones, as in \cite{FernandesFS:2018}. Ultimately, a bound on the total Milnor number of projective hypersurfaces with isolated singularities determines the degree of $P$ completely (Theorem \ref{thm:top}). Crucially, this step relies only on topological properties:

\begin{customthm}{\ref*{thm:top}}
Let $X=V(P)$ and $Y=V(Q)$ be the zero sets of irreducible homogeneous polynomials on $\C^n$, and assume that $\Sing{X}$ has dimension $\leq2$ (e.g., $n\leq4$). If there is a homeomorphism $\Psi\colon\C^n\to\C^n$ with $\Psi(0)=0$ and $\Psi(X)=Y$, then ${\rm deg}(P)={\rm deg}(Q)$.
\end{customthm}

Theorem \ref{thm:top} was previously established for $\dim\Sing{X}=0$ (derived from the topological invariance of $\chi(P^{-1}(1))=1+(-1)^{n-1}({\rm deg}(P)-1)^n$; see the proof of \cite[Theorem 2.2]{FernandesS:2016} and \cite{FernandesJS:2024} for a non-ambient generalization) and for $\dim\Sing{X}\leq 1$, see the proof of \cite[Theorem 5.5]{Sampaio:2020} (see also the proof of \cite[Theorem 6.6]{FernandesS:2023}). In Theorem \ref{thm:main}, the bi-Lipschitz condition is only used in the initial reduction to tangent cones (see Remark \ref{rem:final}(b)).

This technique breaks down when the singular locus of $\mathbb{P} V(P)$ has dimension $\geq 2$. In such scenarios, the Parusi\'nski--Pragacz formula involves the Euler characteristics of higher-dimensional strata intersected with a generic hypersurface of degree ${\rm deg}(P)$. These terms depend on the degrees and Chern--Schwartz--MacPherson classes of the strata closures, which we currently lack the tools to control (furthermore, the degree of a non-hypersurface cone of dimension $\geq3$ is not a bi-Lipschitz invariant for germs \cite{BirbrairFSV:2020, FernandesJS:2023}). Consequently, $n=5$ represents the first open case of Conjecture \ref{conj:metric} outside the scope of Theorem \ref{thm:main}, as it permits tangent cone components with $3$-dimensional singular loci.

\section{Preliminaries}\label{section:preliminaries}

\subsection{Notation and conventions}

Let $A$ be a subset of $\R^N$ with $0$ in its closure $\overline A$. The \textbf{tangent cone} $C(A,0)$ is the set of all vectors $v$ for which there are a sequence of points $x_i\in A\setminus\{0\}$ tending to $0$ and a sequence of positive numbers $t_i$ such that $x_i/t_i \to v$. For a holomorphic germ $f$, the tangent cone of $V(f)$ is $C(V(f),0)=V(in(f))$ (see, e.g., \cite{Chirka:1989}).

Let $P$ be a homogeneous polynomial on $\C^n$. We write $\mathbb{P} V(P)\subset\C\mathbb{P}^{n-1}$ for the projective hypersurface defined by $P$, $F_P=P^{-1}(1)$ for its global Milnor fibre, and $\pi\colon\C^n\setminus\{0\}\to\C\mathbb{P}^{n-1}$ for the standard projection. If $P$ is not constant, $\mathbb{P} V(P)$ is irreducible (resp. reduced) exactly when $P$ is irreducible (resp. reduced). If $P$ is reduced, then, by the Jacobian criterion (see \cite[Section I.5]{Hartshorne:1977}), the singular points of $\mathbb{P} V(P)$ are the points $\pi(z)$ with $\nabla P(z)=0$ (by Euler's formula, $\nabla P(z)=0$ implies $P(z)=0$). We denote by $\Reg{X}=X\setminus\Sing{X}$ the regular (smooth) part of an analytic set $X$.

\medskip
\noindent{\bf Constructible sets.} A subset of $\C^n$ or of $\C\mathbb{P}^n$ is (Zariski) \textbf{locally closed} if it is of the form $A\setminus B$, where $A$ and $B$ are algebraic sets, and it is (complex algebraic) \textbf{constructible} if it is a finite union of locally closed sets. Every algebraic set is constructible, but not conversely: for instance, $\C^*\subset\C$, the complement $\C\mathbb{P}^{n-1}\setminus\mathbb{P} V(P)$ of a projective hypersurface, and $\{x\neq0\}\cup\{(0,0)\}\subset\C^2$ are constructible and not algebraic. Finite unions, finite intersections and complements of constructible sets are constructible. The closure of a constructible set in the classical (Euclidean) topology coincides with its Zariski closure \cite[Chapter I, \S 10, Corollary 1]{Mumford:1988}; in particular, it is an algebraic set.

\medskip
\noindent{\bf Euler characteristics.} We use rational coefficients for homology and cohomology groups unless we state otherwise; $b_i$ is the $i$-th Betti number, and, for a space $A$ whose Betti numbers are finite and vanish in large degrees, $\chi(A)=\sum_i(-1)^ib_i(A)$ is its (singular) Euler characteristic. A complex constructible set is a real semialgebraic set, so it admits a finite semialgebraic triangulation (see \cite[Section 9.2]{BochnakCR:1998}); in particular, its Betti numbers are finite and $\chi$ is defined. For complex constructible sets $A$, $\chi(A)$ coincides with the Euler characteristic with compact supports $\chi_c(A)=\sum_i(-1)^i\dim H^i_c(A)$ (see \cite[p.~141]{Fulton:1993}, \cite[Section 6.0.6]{Schurmann:2003}, \cite[Section 2]{CappellMS:2008}). Since $\chi_c$ is additive, by the long exact sequence of cohomology with compact supports of the pair $(A,B)$, $B\subset A$ closed (see \cite[Section 2]{CappellMS:2008}), we obtain:
\begin{enumerate}[label=(\roman*)]
\item {\it additivity}: $\chi(A)=\chi(A\setminus B)+\chi(B)$ for constructible sets $B\subset A$;
\item {\it Fubini formula}: $\chi(A)=\int_{B}\chi(\rho^{-1}(b))\,d\chi(b)$ for any algebraic map $\rho\colon A\to B$ between constructible sets, where, for a constructible function $\beta$ with finitely many values, $\int_A\beta\, d\chi=\sum_a a\,\chi(A\cap\beta^{-1}(a))$ (see \cite{Viro:1988}). In particular, if all the fibres of $\rho$ have the same Euler characteristic $c$, then $\chi(A)=c\,\chi(B)$; e.g., $\chi(A)=d\,\chi(B)$ if $\rho$ is an unramified covering of degree $d$.
\end{enumerate}
The word ``complex'' is essential here: for an open interval $I\subset\R$, $\chi(I)=1$ but $\chi_c(I)=-1$.

\medskip
\noindent{\bf Milnor fibres.} Let $X$ be a hypersurface in a complex manifold and let $x\in X$. We write $F_{X,x}$ for the Milnor fibre at $x$ of a reduced local equation of $X$ at $x$. It is a smooth manifold of real dimension $2\dim_\C X$, with the homotopy type of a finite CW complex, and its homotopy type does not depend on the choices made in its definition \cite[\S\S 5--6]{Milnor:1968}; it does not depend on the choice of the reduced equation either, by Theorem \ref{thm:le} below applied to the identity map. (All the hypersurfaces to which we apply these notions in Sections \ref{section:cones}--\ref{section:examples} are algebraic, so Milnor's results for polynomials suffice.) We define $c_X(x)=\chi(F_{X,x})$, the singular Euler characteristic of the Milnor fibre. The local Milnor fibre is not a complex constructible set, but, being an orientable manifold of even real dimension, it also satisfies $\chi_c(F_{X,x})=\chi(F_{X,x})$, by Poincar\'e duality with compact supports (see \cite[Section 3.3]{Hatcher:2002}); so there is no ambiguity. If $x$ is an isolated singular point of $X$, then $\mu(X,x)$ is the Milnor number of $X$ at $x$, and $F_{X,x}$ has the homotopy type of a bouquet of $\mu(X,x)$ spheres of dimension $\dim_\C X$ \cite[Theorems 6.5 and 7.2]{Milnor:1968}, so that $c_X(x)=1+(-1)^{\dim X}\mu(X,x)$.

For any integer $n\geq3$ and real number $x>0$, we define
\begin{equation}\label{eq:phi}
\phi_n(x)=\frac{(1-x)^n-1}{x}.
\end{equation}
Then the Euler characteristic of a smooth hypersurface of degree $d$ in $\C\mathbb{P}^{n-1}$ is $\chi_{n,d}:=\phi_n(d)+n$ (see, e.g., \cite[Chapter 5]{Dimca:1992}).

We will use the following elementary fact several times.

\begin{lemma}\label{lemma:product}
Let $p$ be a polynomial on $\C^a$ with $p(u_0)=0$, and let $q(u,t)=p(u)$ on $\C^a\times\C^b$. Then the Milnor fibre of $q$ at $(u_0,t_0)$ is homotopy equivalent to the Milnor fibre of $p$ at $u_0$.
\end{lemma}
\begin{proof}
Take $\varepsilon>0$ such that the Milnor fibrations of $p$ at $u_0$ and of $q$ at $(u_0,t_0)$ are defined on the open balls of radius $\varepsilon$ \cite[\S 5]{Milnor:1968}, and take $\eta\neq0$ small enough for both of them. Then $F_q=\{(u,t):\|u-u_0\|^2+\|t-t_0\|^2<\varepsilon^2,\ p(u)=\eta\}$. The homotopy $H_\lambda(u,t)=(u,t_0+\lambda(t-t_0))$, $\lambda\in[0,1]$, maps $F_q$ into $F_q$ (it does not increase the distance to $(u_0,t_0)$ and it does not change $u$), $H_1={\rm id}$, $H_\lambda$ is the identity on $F_q\cap\{t=t_0\}$, and $H_0(F_q)=F_q\cap\{t=t_0\}$. Hence, $F_q$ deformation retracts onto $F_q\cap\{t=t_0\}=\{u:\|u-u_0\|<\varepsilon,\ p(u)=\eta\}\times\{t_0\}$, which is the Milnor fibre of $p$ at $u_0$.
\end{proof}

\subsection{Topological invariants}

\begin{theorem}[A'Campo \cite{ACampo:1973}, L\^e \cite{Le:1973}]\label{thm:acampo_le}
Let $(X,x)$ and $(Y,y)$ be germs of (reduced) hypersurfaces in $\C^n$, and let $\Psi\colon(\C^n,x)\to(\C^n,y)$ be a homeomorphism of germs with $\Psi(X)=Y$. Then $X$ is smooth at $x$ if and only if $Y$ is smooth at $y$.
\end{theorem}

\begin{definition}\label{def:coherent}
Let $(X,x)$ be a germ of hypersurface in $\C^n$ and let $z$ be a regular point of $X$ close to $x$. A {\bf meridian} of $X$ at $z$ is a small positively oriented circle around $z$ in a complex line through $z$ transverse to $X$ at $z$. If $B$ is a small ball centred at $z$, then $X\cap B$ is a smooth connected hypersurface of $B$, the pair $(B,X\cap B)$ is homeomorphic to the pair formed by a ball and a linear hyperplane through its centre (by the conic structure theorem \cite[Theorem 2.10]{Milnor:1968}, since the link of a smooth point is an unknotted sphere), $B\setminus X$ is homotopy equivalent to a circle, and $H_1(B\setminus X;\Z)\cong\Z$ is generated by the class of any meridian at a point of $X\cap B$.

Let $\Psi\colon(\C^n,x)\to(\C^n,y)$ be a homeomorphism of germs with $\Psi(X)=Y$. By Theorem \ref{thm:acampo_le}, $\Psi(z)$ is a regular point of $Y$. Choose small balls $B'$ centred at $\Psi(z)$ and $B$ centred at $z$ as above (for $Y$ and $X$, respectively), with $\Psi(B)\subset B'$, and a ball $B''\subset\Psi(B)$ centred at $\Psi(z)$. The inclusions induce an isomorphism $H_1(B''\setminus Y;\Z)\to H_1(B'\setminus Y;\Z)$ which factors through $H_1(\Psi(B)\setminus Y;\Z)\cong H_1(B\setminus X;\Z)\cong\Z$; hence $\Psi_*\colon H_1(B\setminus X;\Z)\to H_1(B'\setminus Y;\Z)$ is an isomorphism $\Z\to\Z$, and $\Psi$ maps a meridian at $z$ to a loop which is homologous in $B'\setminus Y$ to $\varepsilon(z)$ times a meridian at $\Psi(z)$, with $\varepsilon(z)\in\{1,-1\}$. We say that $\Psi$ is {\bf coherent} at $x$ if $\varepsilon(z)$ is the same for all the regular points $z$ of $X$ close to $x$.
\end{definition}

For instance, a biholomorphism is coherent at every point, with $\varepsilon\equiv1$, since it preserves the complex orientation of the meridians. If $(X,x)$ is irreducible, every homeomorphism $\Psi$ as above is coherent at $x$, since $\varepsilon$ is locally constant (see the proof of Lemma \ref{lemma:coherent}) and the regular part of a small representative of an irreducible germ is connected (see \cite[Chapter 1]{Chirka:1989}).

\begin{theorem}[L\^e {\cite[Theorem 3.3 and Remark 3.4]{Le:1973a}}]\label{thm:le}
Let $(X,x)$ and $(Y,y)$ be germs of hypersurfaces in $\C^n$ with reduced equations $f$ and $g$, and let $\Psi\colon(\C^n,x)\to(\C^n,y)$ be a homeomorphism of germs with $\Psi(X)=Y$ which is coherent at $x$. Then the Milnor fibres $F_{X,x}$ and $F_{Y,y}$ are homotopy equivalent, and the map induced in homology by the monodromy of $g$ is conjugate to the map induced by the monodromy of $f$ or by its inverse. In particular, $\chi(F_{X,x})=\chi(F_{Y,y})$.
\end{theorem}

Let us recall the idea of the proof, which explains the role of coherence. Let $B$ be a small ball centred at $x$. The group $H_1(B\setminus X;\Z)$ is generated by meridians at regular points of $X$, and $f_*\colon H_1(B\setminus X;\Z)\to H_1(\C^*;\Z)=\Z$ maps each meridian to $1$, since $f$ is reduced. L\^e proved in \cite{Le:1973a} that $F_{X,x}$ is homotopy equivalent to the infinite cyclic covering of $B\setminus X$ associated with $f_*$, the monodromy corresponding to a generator of the group of deck transformations. If $\Psi$ is coherent at $x$, with sign $\varepsilon$, then $g_*\circ\Psi_*=\varepsilon f_*$ on $H_1(B\setminus X;\Z)$, so $\Psi$ lifts to a homeomorphism between the corresponding infinite cyclic coverings (the dependence on the neighbourhoods is handled in \cite{Le:1973a} by comparing nested balls), which conjugates the generator of the deck group to the generator ($\varepsilon=1$) or to its inverse ($\varepsilon=-1$). See also \cite[Proposition 6.2]{FernandesS:2023} and \cite[Lemma 2.5]{FernandesJS:2025}.

\begin{lemma}\label{lemma:coherent}
Let $X$ and $Y$ be the zero sets of irreducible homogeneous polynomials on $\C^n$, and let $\Psi\colon\C^n\to\C^n$ be a homeomorphism with $\Psi(0)=0$ and $\Psi(X)=Y$. Then $\Psi(\Sing{X})=\Sing{Y}$, and $\Psi$ is coherent at every point of $X$.
\end{lemma}
\begin{proof}
For $x\in X$, the germ of $\Psi$ at $x$ is a homeomorphism of germs $(\C^n,x)\to(\C^n,\Psi(x))$ mapping $X$ onto $Y$, so, by Theorem \ref{thm:acampo_le}, $x\in\Sing{X}$ if and only if $\Psi(x)\in\Sing{Y}$. Hence $\Psi(\Sing{X})=\Sing{Y}$.

Let $\varepsilon\colon\Reg{X}\to\{1,-1\}$ be the function given by Definition \ref{def:coherent}. We claim that $\varepsilon$ is locally constant. Let $z\in\Reg{X}$, and let $B$ and $B'$ be balls centred at $z$ and $\Psi(z)$ as in Definition \ref{def:coherent}. Let $z'\in\Reg{X}\cap B$ be close to $z$, and choose balls $B_1'\subset B'$ centred at $\Psi(z')$ and $B_1\subset B$ centred at $z'$ as in Definition \ref{def:coherent}, with $\Psi(B_1)\subset B_1'$. Let $m_z$, $m_{z'}$, $m_{\Psi(z)}$ and $m_{\Psi(z')}$ be meridians at these points, contained in $B$, $B_1$, $B'$ and $B_1'$, respectively. Since $H_1(B\setminus X;\Z)\cong\Z$ is generated by any meridian at a point of $X\cap B$, $m_z$ and $m_{z'}$ are homologous in $B\setminus X$, and similarly $m_{\Psi(z)}$ and $m_{\Psi(z')}$ are homologous in $B'\setminus Y$. Hence, in $H_1(B'\setminus Y;\Z)$,
$$
\varepsilon(z)[m_{\Psi(z)}]=[\Psi(m_z)]=[\Psi(m_{z'})]=\varepsilon(z')[m_{\Psi(z')}]=\varepsilon(z')[m_{\Psi(z)}],
$$
so $\varepsilon(z)=\varepsilon(z')$. Since $X$ is irreducible, $\Reg{X}$ is connected (see, e.g., \cite[Chapter 1]{Chirka:1989}), hence $\varepsilon$ is constant. Finally, for any $x\in X$, the regular points of $X$ close to $x$ are regular points of $X$, so $\Psi$ is coherent at $x$.
\end{proof}

The following result is essentially known: the equality $d=e$ when $\chi(F_P)\neq0$ was obtained in the proof of \cite[Theorem 2.3]{FernandesS:2016} (see Claim 2 there; the argument is written for $n=3$, but it works for every $n$), see also \cite[Proposition 5.3]{Sampaio:2020} and \cite[Claim 6.6.1 and Case 1 in the proof of Theorem 6.6]{FernandesS:2023}, and the equality $\chi(\mathbb{P} V(P))=\chi(\mathbb{P} V(Q))$ holds in general by \cite{FernandesJS:2025}. Since the proof is short, we present it here.

\begin{proposition}\label{prop:lefschetz}
Let $P$ and $Q$ be irreducible homogeneous polynomials on $\C^n$ of degrees $d$ and $e$, and let $\Psi\colon\C^n\to\C^n$ be a homeomorphism with $\Psi(0)=0$ and $\Psi(V(P))=V(Q)$. Then $\chi(F_P)=\chi(F_Q)$. Moreover:
\begin{enumerate}[label=(\alph*)]
\item if $\chi(F_P)\neq0$, then $d=e$;
\item if $\chi(F_P)=0$, then $\chi(\mathbb{P} V(P))=\chi(\mathbb{P} V(Q))=n$.
\end{enumerate}
\end{proposition}
\begin{proof}
The map $P\colon\C^n\setminus V(P)\to\C^*$ is a locally trivial fibration, whose fibre is $F_P$ and whose geometric monodromy is $h(z)=e^{2\pi i/d}z$; it is equivalent to the local Milnor fibration of $P$ at $0$ (see \cite[\S 9]{Milnor:1968} and \cite[Chapter 3]{Dimca:1992}). The same holds for $Q$, with monodromy $h_Q(z)=e^{2\pi i/e}z$. Let $\Lambda(h^k)=\sum_q(-1)^q\,{\rm tr}\big(h^k_*\colon H_q(F_P)\to H_q(F_P)\big)$ be the Lefschetz number of $h^k$.

Since $h^d={\rm id}$, $\Lambda(h^k)$ depends only on the class of $k$ modulo $d$, and $\Lambda(h^k)=\chi(F_P)$ if $d$ divides $k$. If $d$ does not divide $k$, then $h^k$ has no fixed points on $F_P$ (since $h^k(z)=z$ implies $z=0\notin F_P$), and $\Lambda(h^k)=0$ (see \cite[Claim 2 in the proof of Theorem 2.3]{FernandesS:2016}, where $F_P$ is replaced by the compact $h$-invariant deformation retract $F_P\cap \overline{B}_R$, $R\gg1$, to which the Lefschetz fixed point theorem applies; for $1\leq k<d=m(V(P),0)$ this also follows from A'Campo's formula for the Lefschetz numbers of the iterates of the monodromy \cite{ACampo:1975}). Thus
$$
\Lambda(h^k)=\chi(F_P)\ \mbox{ if } d\mid k,\qquad \Lambda(h^k)=0\ \mbox{ otherwise},
$$
and the same holds for $Q$, with $e$ in place of $d$.

By Lemma \ref{lemma:coherent}, $\Psi$ is coherent at $0$, so, by Theorem \ref{thm:le}, $F_P$ and $F_Q$ are homotopy equivalent, in particular $\chi(F_P)=\chi(F_Q)$, and $h_{Q*}$ is conjugate to $h_*^{\varepsilon}$, $\varepsilon\in\{1,-1\}$. Since $h_*^d={\rm id}$, the eigenvalues of $h_*$ on $H_q(F_P)$ are roots of unity, so ${\rm tr}(h_*^{-k})=\overline{{\rm tr}(h_*^k)}={\rm tr}(h_*^k)$, the trace being rational. Hence $\Lambda(h_Q^k)=\Lambda(h^{\varepsilon k})=\Lambda(h^k)$ for all $k\geq1$.

(a) If $\chi(F_P)\neq0$, then $\chi(F_Q)\neq0$, and the set of positive integers $k$ with $\Lambda(h^k)\neq0$ is $d\Z_{>0}$; it must be equal to $e\Z_{>0}$, hence $d=e$.

(b) The map $F_P\to\C\mathbb{P}^{n-1}\setminus\mathbb{P} V(P)$, $z\mapsto\pi(z)$, is the quotient map of the free action of the cyclic group generated by $h$, so it is an unramified covering of degree $d$ between complex constructible sets. Hence, by the Fubini formula and the additivity of $\chi$, $\chi(F_P)=d\,\chi(\C\mathbb{P}^{n-1}\setminus\mathbb{P} V(P))=d\,(n-\chi(\mathbb{P} V(P)))$, and similarly $\chi(F_Q)=e\,(n-\chi(\mathbb{P} V(Q)))$. If $\chi(F_P)=0$, then $\chi(F_Q)=0$, and so $\chi(\mathbb{P} V(P))=\chi(\mathbb{P} V(Q))=n$.
\end{proof}

\subsection{Bi-Lipschitz invariants}

\begin{theorem}[{\cite[Theorem 3.2]{Sampaio:2016}, \cite[Proposition 1.6]{FernandesS:2016}; see also \cite[Theorems 3.11, 5.1 and 5.3]{FernandesS:2023}}]\label{thm:tangent_cones}
Let $f,g\colon(\C^n,0)\to(\C,0)$ be germs of reduced holomorphic functions, and let $\varphi\colon(\C^n,0)\to(\C^n,0)$ be a bi-Lipschitz homeomorphism such that $\varphi(V(f))=V(g)$. Then there is a bi-Lipschitz homeomorphism $\Psi\colon\C^n\to\C^n$ with $\Psi(0)=0$ and $\Psi(C(V(f),0))=C(V(g),0)$. Moreover, if $X_1,\dots,X_s$ are the irreducible components of $C(V(f),0)$, then $Y_i:=\Psi(X_i)$, $i=1,\dots,s$, are the irreducible components of $C(V(g),0)$, and $ k_{V(f)}(X_i)= k_{V(g)}(Y_i)$ for all $i$.
\end{theorem}

Here $ k_{V(f)}(X_i)$ is the relative multiplicity (Lelong number) of $X_i$ defined by Kurdyka and Raby \cite{KurdykaR:1989} (see \cite{FernandesS:2016}), and we have
\begin{equation}\label{eq:mult_relative}
m(V(f),0)=\sum_{i=1}^s k_{V(f)}(X_i)\,{\rm deg}(X_i)
\end{equation}
(see \cite{KurdykaR:1989}, \cite[Equation (1)]{FernandesS:2016}, \cite[Proposition 5.2]{FernandesS:2023}).

Let us recall how $\Psi$ is obtained. For small $t>0$, the maps $\varphi_t(v)=t^{-1}\varphi(tv)$ are bi-Lipschitz with the same constants as $\varphi$, and they are defined on balls whose radii tend to $+\infty$ as $t\to0^+$. By the Arzel\`a--Ascoli theorem, there is a sequence $t_k\to0^+$ such that $\varphi_{t_k}$ and $(\varphi^{-1})_{t_k}$ converge uniformly on compact sets to maps $\Psi$ and $\Psi'$, respectively, and then $\Psi'=\Psi^{-1}$; so $\Psi$ is bi-Lipschitz and $\Psi(0)=0$. Since the tangent cone of a subanalytic set $A$ at $0$ is the limit of the sets $t^{-1}A$ as $t\to0^+$ (in the Hausdorff sense, on compact sets), $\Psi(C(V(f),0))=C(V(g),0)$ \cite{Sampaio:2016}. By Theorem \ref{thm:acampo_le} applied at every point of the reduced hypersurface $C(V(f),0)$ (or by the Lipschitz regularity theorem \cite{BirbrairFLS:2016, Sampaio:2016}), $\Psi$ maps the singular points of $C(V(f),0)$ onto the singular points of $C(V(g),0)$. Since the irreducible components of a hypersurface are the closures of the connected components of its regular part, $\Psi$ maps the irreducible components of $C(V(f),0)$ onto those of $C(V(g),0)$. The equality of the relative multiplicities is \cite[Proposition 1.6]{FernandesS:2016} (its proof applies to the restriction of $\Psi$ to $C(V(f),0)$, which is the limit of the restrictions of $\varphi_{t_k}$ to $t_k^{-1}V(f)$).

\subsection{Whitney stratifications}\label{subsection:whitney}

Let $M\subset\C\mathbb{P}^{N+1}$ be an algebraic set. By a \textbf{Whitney stratification of $M$ with algebraic strata} we mean a finite partition $\mathcal Z$ of $M$ into connected complex submanifolds of $\C\mathbb{P}^{N+1}$ (the \textbf{strata}) which are Zariski locally closed, such that every pair of strata $(A,S)$ satisfies Whitney's conditions (a) and (b) at every point of $S\cap\overline A$. Such a stratification exists, and it can be chosen so that $\Sing{M}$ is a union of strata (see \cite{Whitney:1965}, \cite{Verdier:1976}; the connected components of a smooth Zariski locally closed set are Zariski locally closed, since they are its irreducible components). We do not use the frontier condition. The strata of $\mathcal Z$, together with $\C\mathbb{P}^{N+1}\setminus M$, form a Whitney stratification $\mathcal Z'$ of $\C\mathbb{P}^{N+1}$: Whitney's conditions for the pairs $(\C\mathbb{P}^{N+1}\setminus M,S)$ are automatic, since the tangent spaces of the open stratum are the whole tangent spaces of $\C\mathbb{P}^{N+1}$.

\begin{lemma}\label{lemma:frontier}
Let $\mathcal Z$ be a Whitney stratification of $M$ with algebraic strata, and let $S\in\mathcal Z$ be a stratum of complex dimension $1$. Then $\overline S$ is an irreducible projective curve, and $\overline S\setminus S$ is a finite set.
\end{lemma}
\begin{proof}
Since $S$ is constructible, its closure $\overline S$ coincides with its Zariski closure \cite[Chapter I, \S 10, Corollary 1]{Mumford:1988}, so $\overline S$ is an algebraic set of dimension $1$, and $S$ is Zariski open in $\overline S$ (write $S=A\setminus B$ with $A,B$ algebraic; then $S=\overline S\setminus B$). Let $C_1,\dots,C_r$ be the irreducible components of $\overline S$. The sets $S\cap C_j$ are closed analytic subsets of the connected complex manifold $S$ which cover $S$. A proper closed analytic subset of a connected complex manifold is nowhere dense (see \cite[Chapter 1]{Chirka:1989}), so a finite union of proper ones is nowhere dense; hence $S\subset C_j$ for some $j$, and then $\overline S\subset C_j$, i.e., $\overline S=C_j$ is irreducible. Finally, $\overline S\setminus S=\overline S\cap B$ is a Zariski closed subset of the irreducible curve $\overline S$ which does not contain $S$; hence it has dimension $0$, i.e., it is finite.
\end{proof}

The following lemma is a particular case of the local topological triviality of Whitney stratifications along strata (see \cite[Part I, Chapter 1]{GoreskyM:1988}, \cite{Mather:2012}, \cite{GibsonWPL:1976}); we state the precise version we need and indicate its proof.

\begin{lemma}\label{lemma:product_structure}
Let $\mathcal Z$ be a Whitney stratification of $M$ with algebraic strata, let $S\in\mathcal Z$ be a stratum of dimension $1$ and $x\in S$. Work in an affine chart $\C^{N+1}$ of $\C\mathbb{P}^{N+1}$ containing $x$, and let $H=\rho^{-1}(0)$ be an affine hyperplane through $x$, where $\rho$ is an affine function, such that $T_xS\not\subset T_xH$. Then there are $\varepsilon,\delta>0$ and a homeomorphism
$$
\Theta\colon (H\cap B_\varepsilon(x))\times D_\delta\to U':=B_\varepsilon(x)\cap\rho^{-1}(D_\delta),
$$
where $D_\delta=\{c\in\C:|c|<\delta\}$, such that:
\begin{enumerate}[label=(\alph*)]
\item $\rho(\Theta(z,c))=c$;
\item $\Theta(z,0)=z$ for all $z\in H\cap B_\varepsilon(x)$;
\item $\Theta\big((A\cap H\cap B_\varepsilon(x))\times D_\delta\big)=A\cap U'$ for every $A\in\mathcal Z'$. In particular, $\Theta\big((M\cap H\cap B_\varepsilon(x))\times D_\delta\big)=M\cap U'$ and $\Theta(\{x\}\times D_\delta)=S\cap U'$.
\end{enumerate}
\end{lemma}
\begin{proof}
{\it Step 1: transversality near $x$.} We claim that there is $\varepsilon_0>0$ such that $H$ is transverse to every stratum of $\mathcal Z'$ at every point of $H\cap \overline B_{\varepsilon_0}(x)$. The strata whose closures do not contain $x$ are at positive distance from $x$, and $x$ belongs to no stratum of dimension $0$. Let $A\in\mathcal Z'$ with $x\in\overline A$, and assume that there are points $y_k\in A\cap H$, $y_k\to x$, at which $H$ is not transverse to $A$, i.e., $T_{y_k}A\subset T_{y_k}H=T_xH$ (as $H$ is a hyperplane). Passing to a subsequence, $T_{y_k}A$ converges to a linear space $\tau\subset T_xH$, and, by Whitney's condition (a) for the pair $(A,S)$ at $x$ (or trivially, if $A=S$), $T_xS\subset\tau\subset T_xH$, a contradiction. In particular, $S\cap H\cap\overline B_{\varepsilon_0}(x)=\{x\}$ if $\varepsilon_0$ is small, as $S$ and $H$ have complementary dimensions and are transverse at $x$.

{\it Step 2: choice of $\varepsilon$.} By Step 1, near $x$ the sets $A\cap H$, $A\in\mathcal Z'$, form a Whitney stratification of $H$, in which $\{x\}=S\cap H$ is a stratum (transverse intersections of Whitney stratifications with submanifolds are Whitney stratifications, see \cite[Part I, Chapter 1]{GoreskyM:1988}). Small spheres centred at a point of a Whitney stratified set are transverse to all the strata (see \cite[Part I, Chapter 1]{GoreskyM:1988}); so we can fix $0<\varepsilon\leq\varepsilon_0$ such that the sphere $\Sp_\varepsilon(x)$ is transverse in $H$ to every stratum $A\cap H$. Then $\Sp_\varepsilon(x)$ is transverse to $A$ at the points of $A\cap H\cap\Sp_\varepsilon(x)$, and, if $r$ denotes the distance to $x$, the restriction of $\rho$ to $A\cap\Sp_\varepsilon(x)$ is a submersion at these points: indeed, if $V=T_yA$, $K=\ker(d\rho|_V)=T_y(A\cap H)$ and $W=\ker(dr|_V)$, then $dr|_K\neq0$ gives $W+K=V$, so $d\rho(W)=d\rho(V)=\C$.

{\it Step 3: a stratified submersion.} Stratify the closed ball $\overline B=\overline B_\varepsilon(x)$ by the sets $A\cap B_\varepsilon(x)$ and $A\cap\Sp_\varepsilon(x)$, $A\in\mathcal Z'$. By Steps 1 and 2, the restriction of $\rho$ to each of these strata is a submersion at every point of $\rho^{-1}(0)\cap\overline B$. The set of points of $\overline B$ at which the restriction of $\rho$ (resp. of $r$) to the stratum containing them is not a submersion is closed, by Whitney's condition (a): if $y_k\to y$, with $y_k$ in a stratum $A$, $y$ in a stratum $A'$, and $d\rho|_{T_{y_k}A}$ not surjective, then $d\rho$ is not surjective on a limit $\tau$ of $T_{y_k}A$, and $T_yA'\subset\tau$. By compactness, there is $\delta>0$ such that the restriction of $\rho$ to every stratum is a submersion at every point of $\overline B\cap\rho^{-1}(\overline D_\delta)$, and $\Sp_\varepsilon(x)$ is still transverse to the strata of $\mathcal Z'$ there. The map $\rho\colon\overline B\cap\rho^{-1}(D_\delta)\to D_\delta$ is proper.

{\it Step 4: Thom's first isotopy lemma.} By Thom's first isotopy lemma (see \cite[\S 11]{Mather:2012}, \cite[Chapter II]{GibsonWPL:1976}), the proper stratified submersion of Step 3 is a topologically trivial fibration over the disc $D_\delta$ (the trivialization is obtained by integrating controlled lifts of the coordinate vector fields of $D_\delta$): there is a homeomorphism $\Theta_0\colon(\rho^{-1}(0)\cap\overline B)\times D_\delta\to\overline B\cap\rho^{-1}(D_\delta)$ with $\rho\circ\Theta_0={\rm pr}_2$, which maps $(\Gamma\cap\rho^{-1}(0))\times D_\delta$ onto $\Gamma\cap\rho^{-1}(D_\delta)$ for every stratum $\Gamma$ of $\overline B$. Put $\Theta(z,c)=\Theta_0\big(\Theta_0(\cdot,0)^{-1}(z),c\big)$; then (a) and (b) hold, and $\Theta$ still preserves the strata. Since the strata contained in $\Sp_\varepsilon(x)$ are preserved, $\Theta$ maps $(H\cap B_\varepsilon(x))\times D_\delta$ onto $B_\varepsilon(x)\cap\rho^{-1}(D_\delta)=U'$, and (c) follows. Finally, $S\cap H\cap B_\varepsilon(x)=\{x\}$, so $\Theta(\{x\}\times D_\delta)=S\cap U'$.
\end{proof}

\begin{corollary}\label{cor:constancy}
Let $\mathcal Z$ be a Whitney stratification of a hypersurface $M\subset\C\mathbb{P}^{N+1}$ with algebraic strata, and let $c_M(y)=\chi(F_{M,y})$. Then $c_M$ is constant on every stratum of dimension $1$ of $\mathcal Z$.
\end{corollary}
\begin{proof}
Let $S$ be a stratum of dimension $1$; since $S$ is connected, it suffices to show that $c_M$ is locally constant on $S$. Let $x\in S$, choose an affine hyperplane $H$ through $x$ with $T_xS\not\subset T_xH$, and let $\Theta$ be as in Lemma \ref{lemma:product_structure}. For $c\in D_{\delta/2}$ and $s\in[0,1]$, let $g_s=\Theta\circ\tau_{sc}\circ\Theta^{-1}$, where $\tau_{sc}(z,c')=(z,c'+sc)$; $g_s$ is a homeomorphism from a neighbourhood of $x$ onto a neighbourhood of $\Theta(x,sc)\in S$ which maps $M$ onto $M$, and $g_0={\rm id}$. Let $z$ be a regular point of $M$ close to $x$ and $m_z$ a meridian at $z$. The loops $g_s(m_z)$ vary continuously with $s$ in the complement of $M$, and the points $g_s(z)$ are regular points of $M$ (Theorem \ref{thm:acampo_le}). If $s$ is close to $s_0$, $g_s(m_z)$ is uniformly close to $g_{s_0}(m_z)$, at a distance smaller than the distance from $g_{s_0}(m_z)$ to $M$, so these loops are homotopic in the complement of $M$ near $g_{s_0}(z)$ (by the straight line homotopy), and the meridians at $g_s(z)$ and $g_{s_0}(z)$ are homologous there. Hence the sign $\varepsilon_s(z)$ of Definition \ref{def:coherent} for $g_s$ at $z$ is locally constant in $s$, so $\varepsilon_1(z)=\varepsilon_0(z)=1$. Thus $g_1$ is coherent at $x$, and, by Theorem \ref{thm:le}, $c_M(\Theta(x,c))=c_M(x)$. Since $\Theta(\{x\}\times D_{\delta/2})$ is a neighbourhood of $x$ in $S$, $c_M$ is locally constant on $S$.
\end{proof}

The same argument, using the local topological triviality along strata of any dimension (see \cite[Part I, Chapter 1]{GoreskyM:1988}, \cite{Mather:2012}), shows that $c_M$ is constant on every stratum of $\mathcal Z$; this is used in the statement of Theorem \ref{thm:PP} below. In Section \ref{section:cones} we only need Corollary \ref{cor:constancy}, since there the strata contained in $\Sing{M}$ have dimension $\leq1$, and $c_M\equiv1$ on $\Reg{M}$ (the Milnor fibre at a smooth point is contractible).

\subsection{The Euler characteristic of singular hypersurfaces}

\begin{theorem}[Parusi\'nski--Pragacz \cite{ParusinskiP:1995}]\label{thm:PP}
Let $M\subset\C\mathbb{P}^{N+1}$ be a reduced hypersurface of degree $d$. Let $G$ be a homogeneous polynomial of degree $d$ such that $\mathbb{P} V(G)\subset\C\mathbb{P}^{N+1}$ is smooth and transverse to all the strata of a Whitney stratification $\mathcal{Z}$ of $M$. Then
$$
\chi(M)-\chi_{N+2,d}=-\sum\limits_{S\in \mathcal{Z}}(\chi(F_S)-1)\,\chi(S\setminus \mathbb{P} V(G)),
$$
where $\chi_{N+2,d}=\frac{(1-d)^{N+2}-1}{d}+N+2$ and $F_S$ is the Milnor fibre of $M$ at any point of $S$.
\end{theorem}

The right-hand side depends only on the function $c_M(y)=\chi(F_{M,y})$, and not on the stratification $\mathcal Z$, which only serves to guarantee that $c_M$ is constructible and to define the genericity of $G$. More precisely:

\begin{corollary}\label{cor:PP_partition}
Under the hypotheses of Theorem \ref{thm:PP}, let $M\setminus\mathbb{P} V(G)=P_1\sqcup\dots\sqcup P_r$ be any finite partition into constructible sets such that $c_M$ is constant on each $P_j$, say equal to $c_j$. Then
$$
\chi(M)-\chi_{N+2,d}=-\int_{M\setminus\mathbb{P} V(G)}(c_M-1)\,d\chi=-\sum_{j=1}^r(c_j-1)\,\chi(P_j).
$$
\end{corollary}
\begin{proof}
The sets $S\setminus\mathbb{P} V(G)$, $S\in\mathcal Z$, also form a finite partition of $M\setminus\mathbb{P} V(G)$ into constructible sets on which $c_M$ is constant, equal to $\chi(F_S)$. Refine both partitions: the sets $P_j\cap(S\setminus\mathbb{P} V(G))$ are constructible, and, if such a set is nonempty, then $c_j=\chi(F_S)$. By the additivity of $\chi$,
\begin{align*}
\sum_{S\in\mathcal Z}(\chi(F_S)-1)\,\chi(S\setminus\mathbb{P} V(G))
&=\sum_{S\in\mathcal Z}\sum_{j=1}^r(\chi(F_S)-1)\,\chi\big(P_j\cap(S\setminus\mathbb{P} V(G))\big)\\
&=\sum_{j=1}^r(c_j-1)\sum_{S\in\mathcal Z}\chi\big(P_j\cap(S\setminus\mathbb{P} V(G))\big)\\
&=\sum_{j=1}^r(c_j-1)\,\chi(P_j),
\end{align*}
and both expressions are equal to $\int_{M\setminus\mathbb{P} V(G)}(c_M-1)\,d\chi$. The result follows from Theorem \ref{thm:PP}.
\end{proof}

If $M$ has only isolated singular points $p_1,\dots,p_k$, then $c_M(p)-1=(-1)^N\mu(M,p)$, a generic $\mathbb{P} V(G)$ does not contain the points $p_i$, and Corollary \ref{cor:PP_partition} gives the classical formula $\chi(M)=\chi_{N+2,d}+(-1)^{N+1}\sum_i\mu(M,p_i)$ (see \cite{Parusinski:1988} and \cite[Chapter 5]{Dimca:1992}).

\begin{lemma}\label{lemma:milnor_bound}
Let $V\subset\C\mathbb{P}^N$, $N\geq1$, be a reduced hypersurface of degree $d$ with finitely many singular points. Then $\sum_{p\in\Sing{V}}\mu(V,p)\leq(d-1)^N$.
\end{lemma}
\begin{proof}
This is well known; it also follows from the formula of Dimca and Papadima for the degree of the gradient map \cite[Section 3]{DimcaP:2003}, and \cite[Theorem 1]{Huh:2014} gives a sharper bound. We include the short proof. If $d=1$, $V$ is smooth, so we may assume $d\geq2$. Since $\Sing{V}$ is finite, and since the hyperplanes tangent to $V$ at some regular point form a set of dimension $\leq N-1$ in the dual projective space (for each regular point there is exactly one tangent hyperplane), a generic hyperplane $H_\infty$ does not meet $\Sing{V}$ and is transverse to $V$, i.e., $V\cap H_\infty$ is a smooth hypersurface of $H_\infty$. Choose homogeneous coordinates $(x_0:\dots:x_N)$ with $H_\infty=\{x_0=0\}$, let $F$ be a reduced equation of $V$ and $g(x)=F(1,x_1,\dots,x_N)$. Then ${\rm deg}(g)=d$, and the homogeneous part of degree $d$ of $g$ is $g_d=F(0,x_1,\dots,x_N)$, which has an isolated critical point at $0\in\C^N$, because $V\cap H_\infty=\{g_d=0\}$ is smooth (by Euler's formula, a common zero $x\neq0$ of the partial derivatives of $g_d$ would be a singular point of $\{g_d=0\}$). Since $\partial g/\partial x_i(x)=\partial F/\partial x_i(1,x)$ and $\partial F/\partial x_i(0,x)=\partial g_d/\partial x_i(x)$, the hypersurfaces $V(\partial F/\partial x_i)\subset\C\mathbb{P}^N$, $i=1,\dots,N$, which have degree $d-1$, have no common point on $H_\infty$. Hence their intersection is a projective algebraic set contained in the affine chart $\C^N=\C\mathbb{P}^N\setminus H_\infty$, so it is finite, and it is the set of critical points of $g$. By B\'ezout's theorem (see, e.g., \cite[Chapters 7 and 8]{Fulton:1998}; the local rings being Cohen--Macaulay, the intersection multiplicities are the lengths of the local quotient rings), the sum over the critical points $q$ of $g$ of $\dim_\C\mathcal O_{\C^N,q}/(\partial g/\partial x_1,\dots,\partial g/\partial x_N)=\mu(g,q)$ is $(d-1)^N$. Since $\Sing{V}$ is contained in the affine chart, each $p\in\Sing{V}$ is a critical point of $g$, and $\mu(V,p)=\mu(g,p)$ ($g$ being a reduced local equation of $V$ at $p$), the result follows.
\end{proof}

\begin{lemma}\label{lemma:derivative}
Let $n\geq3$. Then $|\phi_n'(x)|>(x-1)^{n-2}$ for all real $x>2$.
\end{lemma}
\begin{proof}
We have $\phi_n'(x)=\big(1-(1-x)^{n-1}((n-1)x+1)\big)/x^2$. Let $y=x-1>1$. Then $|\phi_n'(x)|\geq\big(y^{n-1}((n-1)x+1)-1\big)/x^2$, and
\begin{align*}
y^{n-1}((n-1)x+1)-y^{n-2}x^2&=y^{n-2}\big(y((n-1)(y+1)+1)-(y+1)^2\big)\\
&=y^{n-2}\big((n-2)y^2+(n-2)y-1\big)>1,
\end{align*}
since $n\geq3$ and $y>1$. Hence $y^{n-1}((n-1)x+1)-1>y^{n-2}x^2$, i.e., $|\phi_n'(x)|>y^{n-2}$.
\end{proof}

\section{Homogeneous hypersurfaces with low dimensional singular locus}\label{section:cones}

Throughout this section, $P$ is an irreducible homogeneous polynomial of degree $d\geq2$ on $\C^n$, $n\geq3$, $X=V(P)$, $M=\mathbb{P} V(P)\subset\C\mathbb{P}^{n-1}$, and $N=n-2=\dim_\C M$. Since $P$ is irreducible, $P$ is a reduced equation of $X$ at each of its points (if the germ of $P$ at some point had a multiple factor $h^2$, the partial derivatives of $P$ would vanish on the hypersurface germ $V(h)$; but $V(P,\partial P/\partial z_1,\dots,\partial P/\partial z_n)$ is a proper algebraic subset of the irreducible hypersurface $X$, since $P$ does not divide the nonzero polynomials $\partial P/\partial z_j$ of degree $d-1$).

\medskip
\noindent{\bf Local product structure.} Let $x\in X\setminus\{0\}$, say with $x_1\neq0$. The map $\Phi(z)=(z_1,z_2/z_1,\dots,z_n/z_1)$ is a biholomorphism from $\{z_1\neq0\}$ onto $\C^*\times\C^{n-1}$, and $P(z)=z_1^d\,p(z_2/z_1,\dots,z_n/z_1)$, where $p(u)=P(1,u)$ is a reduced equation of $M$ in the affine chart $\{z_1\neq0\}\cong\C^{n-1}$ of $\C\mathbb{P}^{n-1}$. Hence $\Phi(X\cap\{z_1\neq0\})=\C^*\times V(p)$, and $\Phi$ maps the germ $(\C^n,X,x)$ biholomorphically onto $(\C\times\C^{n-1},\C\times V(p),(x_1,\pi(x)))$. Since $\Phi$ is coherent (Definition \ref{def:coherent}), Theorem \ref{thm:le} and Lemma \ref{lemma:product} imply that $F_{X,x}$ is homotopy equivalent to $F_{M,\pi(x)}$. Therefore
\begin{equation}\label{eq:c_cone}
c_X(x)=c_M(\pi(x))\qquad\mbox{and}\qquad \Sing{X}\setminus\{0\}=\pi^{-1}(\Sing{M}),
\end{equation}
and, near $x$, $\Sing{X}$ is biholomorphic to the product of a germ of $\C$ and the germ $(\Sing{M},\pi(x))$.

\medskip
\noindent{\bf The singular curves.} Assume that $\dim\Sing{M}\leq1$ (equivalently, $\dim\Sing{X}\leq2$), let $D_1,\dots,D_k$ be the irreducible components of $\Sing{M}$ of dimension $1$, and let $\delta_i={\rm deg}(D_i)$. Let $\Sigma_0$ be the set of points $y\in\Sing{M}$ such that either $y$ is an isolated point of $\Sing{M}$, or the germ of curve $(\Sing{M},y)$ has at least two branches, or $c_M$ is not constant on a neighbourhood of $y$ in $\Sing{M}$.

The set $\Sigma_0$ is finite. Indeed, fix a Whitney stratification $\mathcal Z$ of $M$ with algebraic strata such that $\Sing{M}$ is a union of strata (Subsection \ref{subsection:whitney}). The strata contained in $\Sing{M}$ have dimension $0$ or $1$. Let $T$ be the union of the strata of dimension $0$ and of the sets $\overline S\setminus S$, where $S$ runs over the strata of dimension $1$ contained in $\Sing{M}$; by Lemma \ref{lemma:frontier}, $T$ is finite. Let $y\in\Sing{M}\setminus T$. Then $y$ belongs to a stratum $S$ of dimension $1$. If $S'\neq S$ is another stratum contained in $\Sing{M}$, then $y\notin\overline{S'}$: if $\dim S'=0$, then $\overline{S'}=S'\subset T$; if $\dim S'=1$, then $\overline{S'}\subset S'\cup T$ and $S\cap S'=\emptyset$. Since $\Sing{M}$ is the finite union of the closed sets $\overline{S'}$, $S'\subset\Sing{M}$, and $\overline S\setminus S\subset T$, there is a neighbourhood $U$ of $y$, disjoint from the finite set $T$ and from the sets $\overline{S'}$, $S'\neq S$, so that $\Sing{M}\cap U=S\cap U$, a smooth connected curve if $U$ is small. Hence $y$ is not isolated in $\Sing{M}$, $(\Sing{M},y)$ has one branch, and $c_M$ is constant near $y$ on $\Sing{M}$, by Corollary \ref{cor:constancy}. So $\Sigma_0\subset T$.

A point $y\in\Sing{M}\setminus\Sigma_0$ belongs to exactly one curve $D_i$ (a point of $D_i\cap D_j$, $i\neq j$, is a point where $\Sing{M}$ has at least two branches), and $(\Sing{M},y)=(D_i,y)$ is an irreducible germ of curve. Let $D_i^\circ=D_i\setminus\Sigma_0$, a Zariski locally closed set. Then
$$
\Sing{M}=\Sigma_0\sqcup D_1^\circ\sqcup\cdots\sqcup D_k^\circ,
$$
and $c_M$ is locally constant on each $D_i^\circ$. By Lemma \ref{lemma:formula}(a) below (whose proof does not use the numbers $c_i$), $D_i^\circ$ is connected, so $c_M$ is constant on $D_i^\circ$, equal to $c_i$, say. Note that this decomposition of $\Sing{M}$ need not be a Whitney stratification (a point of $D_i^\circ$ may be a stratum of dimension $0$ of every Whitney stratification of $M$); this is irrelevant for us, by Corollary \ref{cor:PP_partition}.

\begin{lemma}\label{lemma:formula}
With the above notation:
\begin{enumerate}[label=(\alph*)]
\item $D_i^\circ$ is homeomorphic to the complement of a finite set in a compact connected Riemann surface; in particular, $D_i^\circ$ is connected;
\item $\chi(D_i^\circ)=2-b_1(\pi^{-1}(D_i^\circ))$;
\item we have
\begin{equation}\label{eq:formula}
\chi(M)=\chi_{n,d}-\sum_{i=1}^k(c_i-1)\big(\chi(D_i^\circ)-d\,\delta_i\big)-\sum_{y\in\Sigma_0}(c_M(y)-1).
\end{equation}
\end{enumerate}
\end{lemma}
\begin{proof}
\noindent\textit{Proof of (a).} Let $\nu_i\colon\widetilde D_i\to D_i$ be the normalization of the irreducible projective curve $D_i$. Then $\widetilde D_i$ is a compact connected Riemann surface, $\nu_i$ is a finite surjective holomorphic map, and, for any $y\in D_i$, the number of points of $\nu_i^{-1}(y)$ equals the number of branches of the germ $(D_i,y)$ (see, e.g., \cite[Chapter 1]{Chirka:1989}). A point where $D_i$ has at least two branches is a point where $\Sing{M}$ has at least two branches, so it lies in $\Sigma_0$. Consequently, $\nu_i$ restricts to a continuous bijection from $\widetilde D_i\setminus\nu_i^{-1}(\Sigma_0)$ onto $D_i^\circ$. This restriction is proper, being the restriction of the proper map $\nu_i$ over the open subset $D_i^\circ$ of $D_i$; hence it is a closed map, and thus a homeomorphism. Since $\nu_i^{-1}(\Sigma_0)$ is finite, (a) holds. Notably, $D_i^\circ$ is compact if and only if $D_i\cap\Sigma_0=\emptyset$, in which case $D_i^\circ=D_i$ is homeomorphic to $\widetilde D_i$.

\medskip
\noindent\textit{Proof of (b).} Write $D=D_i^\circ$, and let $g$ be the genus of $\widetilde D_i$ and $r=\#\nu_i^{-1}(\Sigma_0)\geq0$, so that $D$ is homeomorphic to a closed orientable surface of genus $g$ with $r$ points removed, by (a).

{\it Step 1: reduction to a circle bundle.} Let $K=\pi^{-1}(D)\cap\Sp^{2n-1}$. The homotopy $(z,s)\mapsto z/((1-s)+s\|z\|)$, $s\in[0,1]$, preserves $\pi^{-1}(D)$ (since $\pi(\lambda z)=\pi(z)$ for $\lambda>0$), is the identity for $s=0$ and on $K$, and ends in $K$; so $K$ is a deformation retract of $\pi^{-1}(D)$, and $b_1(\pi^{-1}(D))=b_1(K)$. The map $\pi|_K\colon K\to D$ is the restriction to $D$ of the Hopf fibration $\Sp^{2n-1}\to\C\mathbb{P}^{n-1}$, which is the unit circle bundle of the tautological line bundle $\mathcal O(-1)$. Hence $K\to D$ is the unit circle bundle of the complex line bundle $\mathcal O(-1)|_D$; in particular, it is an oriented circle bundle (the fibres are oriented by the complex structure), and its Euler class is $e=c_1(\mathcal O(-1)|_D)=-\eta|_D$, where $\eta\in H^2(\C\mathbb{P}^{n-1};\Z)$ is the hyperplane class: the Euler class of the underlying oriented real plane bundle of a complex line bundle is its first Chern class, and Euler classes are natural under pull-backs (see \cite[\S\S 9 and 14]{MilnorS:1974}).

{\it Step 2: the Gysin sequence.} The Gysin sequence of the oriented circle bundle $K\to D$ (see \cite[\S 12]{MilnorS:1974} or \cite[Section 4.D]{Hatcher:2002}), with rational coefficients, starts as
$$
0\to H^1(D)\xrightarrow{\ \pi^*\ }H^1(K)\to H^0(D)\xrightarrow{\ \cup e\ }H^2(D).
$$
Since $D$ is connected, $H^0(D)=\Q\cdot1$ and $\cup e$ maps $1$ to the rational Euler class $e_\Q$. Hence
\begin{equation}\label{eq:gysin}
b_1(K)=b_1(D)+1\ \mbox{ if } e_\Q=0,\qquad b_1(K)=b_1(D)\ \mbox{ if } e_\Q\neq0 .
\end{equation}
(The Betti numbers of $D$ and $K$ are finite, since these spaces are homotopy equivalent to finite CW complexes.)

{\it Step 3: $D$ non-compact ($r\geq1$).} A connected non-compact $m$-manifold $W$ satisfies $H_j(W)=0$ for $j\geq m$ (see \cite[Proposition 3.29]{Hatcher:2002}); so $H_2(D)=0$, and $H^2(D)=0$ by the universal coefficient theorem. Hence $e_\Q=0$, and \eqref{eq:gysin} gives $b_1(K)=b_1(D)+1$. Since $b_0(D)=1$ and $b_j(D)=0$ for $j\geq2$, we get $\chi(D)=1-b_1(D)=2-b_1(K)$. 

{\it Step 4: $D$ compact ($r=0$).} Then $D=D_i$ is homeomorphic to $\widetilde D_i$ via $\nu_i$, a closed orientable surface of genus $g$, with fundamental class $[D]=\nu_{i*}[\widetilde D_i]$, and evaluation on $[D]$ is an isomorphism $H^2(D)\to\Q$. We have
$$
\langle e,[D]\rangle=-\langle\eta,\nu_{i*}[\widetilde D_i]\rangle=-\langle\nu_i^*\eta,[\widetilde D_i]\rangle=-\delta_i\neq0 .
$$
Indeed, $\langle\nu_i^*\eta,[\widetilde D_i]\rangle=\deg(\nu_i^*\mathcal O(1))$ is the number of zeros, counted with multiplicities, of the holomorphic section of $\nu_i^*\mathcal O(1)$ induced by a linear form $\ell$ (see \cite[Chapter 1, Section 1]{GriffithsH:1978}); if the hyperplane $\{\ell=0\}$ is generic, these zeros are the $\delta_i={\rm deg}(D_i)$ points of $\nu_i^{-1}(D_i\cap\{\ell=0\})$, all of them simple (see \cite[Lecture 18]{Harris:1992}). (For (b), it suffices to know that $\langle\eta,[D]\rangle>0$, which also follows from the fact that $\eta$ is represented by the Fubini--Study form, whose integral over a complex curve is its area \cite[Chapter 0]{GriffithsH:1978}.) Hence $e_\Q\neq0$, and \eqref{eq:gysin} gives $b_1(K)=b_1(D)=2g$, so $\chi(D)=2-2g=2-b_1(K)$.

In both cases, $\chi(D_i^\circ)=2-b_1(K)=2-b_1(\pi^{-1}(D_i^\circ))$, which proves (b).

\medskip
\noindent\textit{Proof of (c).} Choose a homogeneous polynomial $G$ of degree $d$ such that:
\begin{enumerate}[label=(\roman*)]
\item $\mathbb{P} V(G)$ is smooth and transverse to all the strata of the Whitney stratification $\mathcal Z$ fixed above;
\item $\mathbb{P} V(G)\cap\Sigma_0=\emptyset$;
\item $\mathbb{P} V(G)$ meets each $D_i$ transversally at $d\,\delta_i$ distinct smooth points of $D_i$.
\end{enumerate}
Such a $G$ exists: under the Veronese embedding $v_d$ of degree $d$, the hypersurfaces of degree $d$ correspond to the hyperplanes of the target projective space, $v_d(D_i)$ is an irreducible curve of degree $d\,\delta_i$, and each of the conditions (i)--(iii) holds for the hyperplanes outside a proper Zariski closed subset of the dual projective space (by Bertini's theorem \cite[Theorem II.8.18]{Hartshorne:1977} and Kleiman's transversality theorem \cite{Kleiman:1974} for (i), trivially for (ii), and by \cite[Lecture 18]{Harris:1992} for (iii)). By (ii) and (iii), all the points of $\mathbb{P} V(G)\cap D_i$ lie in $D_i^\circ$.

By Theorem \ref{thm:PP} and Corollary \ref{cor:PP_partition}, applied to $M\subset\C\mathbb{P}^{n-1}$ (so that $N+2=n$) with the partition
$$
M\setminus\mathbb{P} V(G)=\big(\Reg{M}\setminus\mathbb{P} V(G)\big)\sqcup\big(\Sigma_0\setminus\mathbb{P} V(G)\big)\sqcup\bigsqcup_{i=1}^k\big(D_i^\circ\setminus\mathbb{P} V(G)\big),
$$
whose pieces are constructible and on which $c_M$ is constant, equal to $1$, to $c_M(y)$ on each point $y$, and to $c_i$, respectively, we obtain
$$
\chi(M)-\chi_{n,d}=-\sum_{i=1}^k(c_i-1)\,\chi\big(D_i^\circ\setminus \mathbb{P} V(G)\big)-\sum_{y\in\Sigma_0\setminus \mathbb{P} V(G)}(c_M(y)-1).
$$
By (ii), $\Sigma_0\setminus\mathbb{P} V(G)=\Sigma_0$. By (ii) and (iii), $D_i^\circ\cap\mathbb{P} V(G)$ consists of $d\,\delta_i$ points, so, by additivity, $\chi(D_i^\circ\setminus \mathbb{P} V(G))=\chi(D_i^\circ)-d\,\delta_i$. This gives \eqref{eq:formula}.
\end{proof}

\begin{lemma}\label{lemma:transversal}
With the above notation, assume $k\geq1$. Then $N\geq2$, and there are positive integers $\mu_1,\dots,\mu_k$ such that $c_i-1=(-1)^{N-1}\mu_i$ for all $i$ and
$$
\sum_{i=1}^k\mu_i\,\delta_i\leq(d-1)^N.
$$
\end{lemma}
\begin{proof}
Since $M$ is reduced of dimension $N$, its singular locus has dimension $\leq N-1$; as $\dim\Sing{M}=1$, $N\geq2$.

\medskip
\noindent{\it Step 0: choice of a generic hyperplane.} Let $\mathcal Z$ and $T\supset\Sigma_0$ be as above, and let $H\subset\C\mathbb{P}^{N+1}$ be a hyperplane such that:
\begin{enumerate}[label=(\roman*)]
\item $H\cap T=\emptyset$;
\item $H$ meets each $D_i$ transversally at $\delta_i={\rm deg}(D_i)$ distinct smooth points of $D_i$;
\item $H$ is transverse to the complex manifold $\Reg{M}$ and to every stratum of $\mathcal Z$.
\end{enumerate}
Each of these conditions holds for $H$ outside a proper Zariski closed subset of the dual projective space $(\C\mathbb{P}^{N+1})^\vee$, of dimension $N+1$, so a generic $H$ satisfies all of them: (i) because the hyperplanes through a given point form a hyperplane of $(\C\mathbb{P}^{N+1})^\vee$; (ii) by the definition of the degree of a curve (see \cite[Lecture 18]{Harris:1992}); (iii) because, for a smooth locally closed subvariety $Z$ of $\C\mathbb{P}^{N+1}$, the hyperplanes tangent to $Z$ at some point form a set of dimension $\leq N$ (it is the image of the conormal variety of $Z$, which has dimension $N$; see \cite[Lecture 15]{Harris:1992}, or use Kleiman's transversality theorem \cite{Kleiman:1974}). By (i), the points of $D_i\cap H$ do not belong to $T$, hence they belong to strata of dimension $1$ contained in $D_i^\circ$.

Let $L\subset\C^n$ be the linear hyperplane with $H=\mathbb{P}(L)$, and choose coordinates such that $L=\{z_n=0\}$; then $H\cong\C\mathbb{P}^N$ with homogeneous coordinates $(z_1:\dots:z_{n-1})$, and $M\cap H=\mathbb{P} V(p)$, where $p=P|_L=P(z_1,\dots,z_{n-1},0)$. For $\hat y\in L$,
\begin{equation}\label{eq:restriction}
dp(\hat y)=dP(\hat y)|_L .
\end{equation}
Condition (iii) says that, for every $y=\pi(\hat y)\in\Reg{M}\cap H$, the nonzero vector $\nabla P(\hat y)$ is not proportional to the linear form $z_n$ defining $L$, i.e., $dP(\hat y)$ does not vanish on $L$; by \eqref{eq:restriction},
\begin{equation}\label{eq:transv}
y=\pi(\hat y)\in\Reg{M}\cap H\ \Longrightarrow\ dp(\hat y)\neq0.
\end{equation}

\medskip
\noindent{\it Step 1: the singular set of the equation $p$ is finite.} Let $Z=\{\pi(\hat y)\in H: p(\hat y)=0,\ dp(\hat y)=0\}$. By \eqref{eq:transv}, $Z\subset\Sing{M}\cap H$. Since $\Sing{M}=\Sigma_0\sqcup D_1^\circ\sqcup\dots\sqcup D_k^\circ$, by (i) and (ii), $\Sing{M}\cap H=\bigsqcup_i(D_i\cap H)$ has $\delta_1+\dots+\delta_k$ points. Hence $Z$ is finite.

\medskip
\noindent{\it Step 2: $M\cap H$ is a reduced hypersurface of degree $d$ of $H\cong\C\mathbb{P}^N$.} We have $p\neq0$, since an irreducible $P$ of degree $d\geq2$ is not divisible by the linear form $z_n$. If $p$ had a multiple factor, $p=a^2b$ with $a$ non-constant, then $p$ and its partial derivatives $\partial_jp=2a(\partial_ja)b+a^2\partial_jb$ would vanish on $\mathbb{P} V(a)$, which has dimension $N-1\geq1$, so $\mathbb{P} V(a)\subset Z$, contradicting Step 1. Hence $p$ is reduced, and $M\cap H=\mathbb{P} V(p)$ is a reduced hypersurface of degree $d$ of $H$.

\medskip
\noindent{\it Step 3: $\Sing{M\cap H}=\Sing{M}\cap H$.} Since $p$ is reduced, $\Sing{M\cap H}=Z$ by the Jacobian criterion, so $\Sing{M\cap H}\subset\Sing{M}\cap H$ by Step 1. Conversely, if $y=\pi(\hat y)\in\Sing{M}\cap H$, then $dP(\hat y)=0$ ($P$ being reduced), hence $dp(\hat y)=0$ by \eqref{eq:restriction}, and $p(\hat y)=0$, so $y\in Z$. Therefore
$$
\Sing{M\cap H}=\Sing{M}\cap H=\bigsqcup_{i=1}^k(D_i\cap H),
$$
a finite set, with $\#(D_i\cap H)=\delta_i$. In particular, each point of $D_i\cap H$ is an isolated singular point of $M\cap H$.

\medskip
\noindent{\it Step 4: $F_{M,x}\simeq F_{M\cap H,x}$ for $x\in D_i\cap H$.} Let $x\in D_i\cap H$, and let $S$ be the stratum of dimension $1$ containing $x$. Work in an affine chart $\C^{N+1}$ containing $x$, and let $\rho$ be an affine function with $H\cap\C^{N+1}=\rho^{-1}(0)$. By (iii), $T_xS\not\subset T_xH$, so Lemma \ref{lemma:product_structure} gives a homeomorphism $\Theta\colon(H\cap B_\varepsilon(x))\times D_\delta\to U'$ with $\Theta(z,0)=z$ and $\Theta((M\cap H\cap B_\varepsilon(x))\times D_\delta)=M\cap U'$. Identifying $H$ with $\C^N$ by affine coordinates $u$ centred at $x$, and letting $p_x$ be the affine equation of $M\cap H$ in $H$ near $x$ (a reduced equation, by Step 2), $\Theta^{-1}$ induces a homeomorphism of germs
$$
\Theta^{-1}\colon(\C^{N+1},M,x)\to(\C^N\times\C,V(p_x)\times\C,(0,0)),
$$
and $V(p_x)\times\C$ has the reduced equation $q(u,t)=p_x(u)$.

We claim that $\Theta^{-1}$ is coherent at $x$, with sign $+1$. (a) Let $z\in M\cap H$ be close to $x$, $z\neq x$. By Step 3, $z\in\Reg{M}$; by (iii), $H$ is transverse to $M$ at $z$, so $M\cap H$ is smooth at $z$ and $T_z(M\cap H)=T_zM\cap T_zH$. A complex line $\ell\subset H$ through $z$ with $\ell\not\subset T_z(M\cap H)$ is not contained in $T_zM$, so a small positively oriented circle $m_z\subset\ell$ around $z$ is a meridian of $M$ at $z$. Since $\Theta(\cdot,0)$ is the identity, $\Theta^{-1}(m_z)=m_z\times\{0\}$, which is a positively oriented meridian of $V(p_x)\times\C$ at $(z,0)$. So the sign of Definition \ref{def:coherent} at $z$ is $+1$. (b) The sign is locally constant on $\Reg{M}$ near $x$ (see the proof of Lemma \ref{lemma:coherent}, which only uses that the map is a homeomorphism of germs preserving the hypersurfaces). (c) Let $M_1,\dots,M_s$ be the irreducible components of the germ $(M,x)$, in a small ball $B'$. Then $\Reg{M}\cap B'=\bigsqcup_jR_j$, with $R_j=M_j\setminus(\Sing{M_j}\cup\bigcup_{l\neq j}M_l)$, and each $R_j$ is connected (see \cite[Chapter 1]{Chirka:1989}). Each $M_j$ has dimension $N$, so the germ $M_j\cap H$ has dimension $\geq N-1\geq1$ and contains points $z\neq x$ arbitrarily close to $x$; such points are regular points of $M$ (Step 3), hence they lie in $R_j\cap H$. By (a) and (b), the sign is $+1$ on each $R_j$, which proves the claim.

By Theorem \ref{thm:le}, $F_{M,x}$ is homotopy equivalent to the Milnor fibre of $q$ at $(0,0)$, which is homotopy equivalent to $F_{M\cap H,x}$ by Lemma \ref{lemma:product}.

\medskip
\noindent{\it Step 5: the Milnor numbers $\mu_i$.} Let $x\in D_i\cap H$. By Step 3, $x$ is an isolated singular point of the hypersurface $M\cap H$ of $H$, which has dimension $N-1$, so $F_{M\cap H,x}$ has the homotopy type of a bouquet of $\mu(M\cap H,x)\geq1$ spheres of dimension $N-1$ \cite[Theorems 6.5 and 7.2]{Milnor:1968} ($\mu(M\cap H,x)\geq1$ because $x$ is a singular point). Since $x\in D_i^\circ$, Step 4 gives
$$
c_i=c_M(x)=\chi(F_{M\cap H,x})=1+(-1)^{N-1}\mu(M\cap H,x).
$$
Hence $\mu_i:=\mu(M\cap H,x)=(-1)^{N-1}(c_i-1)$ depends only on $i$, and $\mu_i\geq1$.

\medskip
\noindent{\it Step 6: the bound.} By Steps 3 and 5, regrouping the singular points of $M\cap H$ curve by curve,
$$
\sum_{p\in\Sing{M\cap H}}\mu(M\cap H,p)=\sum_{i=1}^k\ \sum_{x\in D_i\cap H}\mu(M\cap H,x)=\sum_{i=1}^k\mu_i\,\delta_i .
$$
By Step 2 and Step 3, $M\cap H$ is a reduced hypersurface of degree $d$ of $H\cong\C\mathbb{P}^N$ with finitely many singular points, so Lemma \ref{lemma:milnor_bound} gives $\sum_{i=1}^k\mu_i\,\delta_i\leq(d-1)^N$.
\end{proof}

\begin{lemma}\label{lemma:matching}
Let $Q$ be an irreducible homogeneous polynomial of degree $e\geq2$ on $\C^n$, $Y=V(Q)$, $M'=\mathbb{P} V(Q)$, and let $\Psi\colon\C^n\to\C^n$ be a homeomorphism with $\Psi(0)=0$ and $\Psi(X)=Y$. Assume that $\dim\Sing{M}\leq1$, and use the above notation for $X$ and the same notation with primes for $Y$. Then $\dim\Sing{M'}\leq1$, $k'=k$, and, after renumbering the curves $D_i'$:
\begin{enumerate}[label=(\alph*)]
\item $\Psi(\pi^{-1}(D_i)\cup\{0\})=\pi^{-1}(D_i')\cup\{0\}$, $c_i=c_i'$ and $\chi(D_i^\circ)=\chi(D_i'^\circ)$ for all $i$;
\item there is a bijection $\Sigma_0\to\Sigma_0'$, $y\mapsto y'$, with $c_M(y)=c_{M'}(y')$;
\item $\delta_i=\delta_i'$ for all $i$.
\end{enumerate}
\end{lemma}
\begin{proof}
\noindent\textit{Step 1.} By Lemma \ref{lemma:coherent}, $\Psi(\Sing{X})=\Sing{Y}$ and $\Psi$ is coherent at every point of $X$. Consequently, by Theorem \ref{thm:le}, $c_Y(\Psi(x))=c_X(x)$ for all $x\in X\setminus\{0\}$ (note that $\Psi(x)\neq0$). Since $\Psi$ is a homeomorphism, the topological dimensions of $\Sing{Y}$ and $\Sing{X}$ coincide, which yields $\dim\Sing{Y}\leq 2$ and, via \eqref{eq:c_cone}, $\dim\Sing{M'}\leq 1$.

\medskip
\noindent\textit{Step 2: a topological characterization of $\widehat\Sigma_0:=\pi^{-1}(\Sigma_0)$.} Let $x\in\Sing{X}\setminus\{0\}$ and $y=\pi(x)$. If $y$ is an isolated point of $\Sing{M}$, the germ $(\Sing{M},y)$ is a point; we set $L_y=\emptyset$ and $r_y=0$. Otherwise, $(\Sing{M},y)$ is a germ of curve with $r_y\geq 1$ branches, and its normalization (see the proof of Lemma \ref{lemma:formula}(a)) shows that it is homeomorphic to the union of $r_y$ discs glued at their centres, i.e., to the cone over $L_y$, a disjoint union of $r_y$ circles. By \eqref{eq:c_cone}, the germ of $\Sing{X}$ at $x$ is homeomorphic to $(\C\times{\rm Cone}(L_y),(0,o))$, where $o$ is the vertex. By the K\"unneth formula for pairs (see \cite[Section 3.B]{Hatcher:2002}) and the exact sequence of the pair $({\rm Cone}(L_y),{\rm Cone}(L_y)\setminus\{o\})$, whose first term is contractible and whose second term is homotopy equivalent to $L_y$,
$$
H_j(\Sing{X},\Sing{X}\setminus\{x\})\cong H_{j-2}({\rm Cone}(L_y),{\rm Cone}(L_y)\setminus\{o\})\cong\widetilde H_{j-3}(L_y)
$$
(with the convention $\widetilde H_{-1}(\emptyset)=\Q$). These are the local homology groups of a point in $\R^4$ (i.e., $\Q$ for $j=4$ and $0$ otherwise) if and only if $L_y$ is a single circle, i.e., if and only if $y$ is not isolated in $\Sing{M}$ and $r_y=1$. Furthermore, $\pi\colon\Sing{X}\setminus\{0\}\to\Sing{M}$ is an open map (a $\C^*$-bundle), and $c_X=c_M\circ\pi$ by \eqref{eq:c_cone}; hence $c_M$ is constant near $y$ in $\Sing{M}$ if and only if $c_X$ is constant near $x$ in $\Sing{X}$. Therefore, $x\in\widehat\Sigma_0$ if and only if either the local homology of $\Sing{X}$ at $x$ is not that of a point of $\R^4$, or $c_X$ is not constant on a neighbourhood of $x$ in $\Sing{X}$. The same holds for $Y$. Since $\Psi$ restricts to a homeomorphism from $\Sing{X}\setminus\{0\}$ onto $\Sing{Y}\setminus\{0\}$, local homology groups are invariant under homeomorphisms, and $c_Y\circ\Psi=c_X$, we get $\Psi(\widehat\Sigma_0)=\widehat\Sigma_0'$.

\medskip
\noindent\textit{Step 3: proof of (a) and (b).} The connected components of $\widehat\Sigma_0$ are the finitely many disjoint closed subsets $\pi^{-1}(y)\cong\C^*$ of $\C^n\setminus\{0\}$, $y\in\Sigma_0$, and the same holds for $Y$. Since $\Psi$ maps connected components onto connected components, it induces a bijection $\Sigma_0\to\Sigma_0'$, $y\mapsto y'$, with $\Psi(\pi^{-1}(y))=\pi^{-1}(y')$, and, for $x\in\pi^{-1}(y)$, $c_M(y)=c_X(x)=c_Y(\Psi(x))=c_{M'}(y')$. This proves (b).

Moreover, $\Sing{X}\setminus(\widehat\Sigma_0\cup\{0\})=\pi^{-1}(\Sing{M}\setminus\Sigma_0)$ is the disjoint union of the sets $\pi^{-1}(D_i^\circ)$. They are connected, by Lemma \ref{lemma:formula}(a), since $\pi$ is a $\C^*$-bundle, and they are closed in this subspace, since $D_i^\circ=D_i\cap(\Sing{M}\setminus\Sigma_0)$; being finitely many, they are also open in it, so they are its connected components. The same holds for $Y$. Since $\Psi$ maps $\Sing{X}\setminus(\widehat\Sigma_0\cup\{0\})$ homeomorphically onto $\Sing{Y}\setminus(\widehat\Sigma_0'\cup\{0\})$, we get $k=k'$ and, after renumbering, $\Psi(\pi^{-1}(D_i^\circ))=\pi^{-1}(D_i'^\circ)$. Since $D_i^\circ$ is dense in $D_i$, taking closures in $\C^n$ gives $\Psi(\pi^{-1}(D_i)\cup\{0\})=\pi^{-1}(D_i')\cup\{0\}$. The identity $c_Y\circ\Psi=c_X$ gives $c_i=c_i'$, and Lemma \ref{lemma:formula}(b) gives $\chi(D_i^\circ)=2-b_1(\pi^{-1}(D_i^\circ))=2-b_1(\pi^{-1}(D_i'^\circ))=\chi(D_i'^\circ)$. This proves (a).

\medskip
\noindent\textit{Step 4: proof of (c).} By (a) and since $\Psi(0)=0$, $\Psi$ maps $\pi^{-1}(D_i)$ homeomorphically onto $\pi^{-1}(D_i')$. The sets $\pi^{-1}(D_i)\cup\{0\}$ and $\pi^{-1}(D_i')\cup\{0\}$ are irreducible homogeneous algebraic surfaces of degrees $\delta_i$ and $\delta_i'$, so (c) follows from \cite[Proposition 3.5]{FernandesFS:2018}, which states that two irreducible homogeneous affine algebraic surfaces whose complements of the origin are homeomorphic have the same degree. For completeness, we present the proof, adapted to our notation. It suffices to show that the torsion subgroup of $H^2(\pi^{-1}(D_i);\Z)$ is $\Z_{\delta_i}$. The space $\pi^{-1}(D_i)$ deformation retracts onto $K_i=\pi^{-1}(D_i)\cap\Sp^{2n-1}$, which is a circle bundle over $D_i$ with Euler class $e=-\eta|_{D_i}$ (as in Step 1 of the proof of Lemma \ref{lemma:formula}(b)).

Let us compute the integral cohomology of $D_i$. Let $\nu_i\colon\widetilde D_i\to D_i$ be the normalization, $B_i\subset D_i$ the finite set of points at which $D_i$ has at least two branches, and $A_i=\nu_i^{-1}(B_i)$. Since $\nu_i$ is a closed continuous surjection which is injective outside $A_i$, $D_i$ is the quotient of the closed orientable surface $\widetilde D_i$, of genus $g$ say, obtained by identifying the points of $A_i$ having the same image. Identifying two $0$-cells of a connected CW complex gives a space homotopy equivalent to its wedge with a circle (see \cite[Chapter 0]{Hatcher:2002}); so $D_i$ is homotopy equivalent to $\widetilde D_i\vee\bigvee_{j=1}^{m_i}\Sp^1$, with $m_i=|A_i|-|B_i|$. Hence $H_1(D_i;\Z)\cong\Z^{2g+m_i}$ is free, $H_2(D_i;\Z)\cong\Z$ is generated by the fundamental class $[D_i]=\nu_{i*}[\widetilde D_i]$, and $H_3(D_i;\Z)=0$. By the universal coefficient theorem (see \cite[Section 3.1]{Hatcher:2002}), $H^1(D_i;\Z)$ is free, $H^3(D_i;\Z)=0$, and evaluation on $[D_i]$ is an isomorphism $H^2(D_i;\Z)\to\Z$. Since $\langle e,[D_i]\rangle=-\langle\nu_i^*\eta,[\widetilde D_i]\rangle=-\delta_i$ (see Step 4 of the proof of Lemma \ref{lemma:formula}(b)), $e$ is $-\delta_i$ times a generator of $H^2(D_i;\Z)$.

The Gysin sequence of $K_i\to D_i$ with integer coefficients (see \cite[\S 12]{MilnorS:1974} or \cite[Section 4.D]{Hatcher:2002}) gives an exact sequence
$$
H^0(D_i;\Z)\xrightarrow{\ \cup e\ }H^2(D_i;\Z)\to H^2(K_i;\Z)\to H^1(D_i;\Z)\to H^3(D_i;\Z)=0.
$$
The first map is multiplication by $-\delta_i$ on $\Z$, so $0\to\Z_{\delta_i}\to H^2(K_i;\Z)\to H^1(D_i;\Z)\to 0$ is exact. Since $H^1(D_i;\Z)$ is free, this sequence splits, and the torsion subgroup of $H^2(\pi^{-1}(D_i);\Z)\cong H^2(K_i;\Z)$ is $\Z_{\delta_i}$. The same holds for $D_i'$, so $\Z_{\delta_i}\cong\Z_{\delta_i'}$, i.e., $\delta_i=\delta_i'$.
\end{proof}

\begin{theorem}\label{thm:top}
Let $P$ and $Q$ be irreducible homogeneous polynomials on $\C^n$ of degrees $d$ and $e$, let $X=V(P)$ and $Y=V(Q)$, and let $\Psi\colon\C^n\to\C^n$ be a homeomorphism with $\Psi(0)=0$ and $\Psi(X)=Y$. If $\dim\Sing{X}\leq2$, then $d=e$. In particular, this holds for any irreducible homogeneous polynomials $P$ and $Q$ on $\C^n$, $n\leq4$.
\end{theorem}
\begin{proof}
By Theorem \ref{thm:acampo_le}, $X$ is smooth at $0$, i.e., $d=1$, if and only if $Y$ is smooth at $0$, i.e., $e=1$. If $n\leq2$, an irreducible homogeneous polynomial has degree $1$. So we may assume $n\geq3$ and $d,e\geq2$. By Lemma \ref{lemma:coherent}, $\Sing{Y}=\Psi(\Sing{X})$, so $\dim\Sing{Y}=\dim\Sing{X}\leq2$. Hence the hypotheses are symmetric in $X$ and $Y$ (replace $\Psi$ by $\Psi^{-1}$), and we may assume $d\leq e$. Note that $\dim\Sing{M}\leq1$ by \eqref{eq:c_cone}.

Assume, by contradiction, that $d<e$. By Proposition \ref{prop:lefschetz}(a), $\chi(F_P)=0$, and then, by Proposition \ref{prop:lefschetz}(b), $\chi(M)=\chi(M')=n$. By Lemma \ref{lemma:matching}, we can apply Lemma \ref{lemma:formula} to $X$ and to $Y$ with the same $k$, $c_i$, $\chi(D_i^\circ)$ and $\sum_{y\in\Sigma_0}(c_M(y)-1)=\sum_{y'\in\Sigma_0'}(c_{M'}(y')-1)$. Let $\delta_i'$ be the degrees of the curves $D_i'$. Subtracting the two formulas \eqref{eq:formula}, and since $\chi_{n,d}-\chi_{n,e}=\phi_n(d)-\phi_n(e)$, we get
$$
0=\chi(M)-\chi(M')=\phi_n(d)-\phi_n(e)+\sum_{i=1}^k(c_i-1)(d\,\delta_i-e\,\delta_i').
$$
By Lemma \ref{lemma:matching}(c), $\delta_i=\delta_i'$, so $\phi_n(e)-\phi_n(d)=(d-e)A$, where $A=\sum_i(c_i-1)\delta_i$. By the mean value theorem, there is $\xi\in(d,e)$ with $\phi_n(e)-\phi_n(d)=\phi_n'(\xi)(e-d)$, so $A=-\phi_n'(\xi)$. Since $\xi>d\geq2$, Lemma \ref{lemma:derivative} gives
$$
|A|=|\phi_n'(\xi)|>(\xi-1)^{n-2}>(d-1)^{n-2}=(d-1)^N.
$$
If $k=0$, then $A=0$, which is a contradiction. If $k\geq1$, then, by Lemma \ref{lemma:transversal}, all the numbers $c_i-1=(-1)^{N-1}\mu_i$ have the same sign, so $|A|=\sum_i\mu_i\delta_i\leq(d-1)^N$, which is also a contradiction. Hence $d=e$.

Finally, if $n\leq4$, then $\Sing{X}$ is a proper algebraic subset of the hypersurface $X$, so $\dim\Sing{X}\leq n-2\leq2$.
\end{proof}

\section{Proof of the main results}\label{section:main}

\begin{theorem}\label{thm:main}
Let $f,g\colon(\C^n,0)\to(\C,0)$ be germs of reduced holomorphic functions, and let $\varphi\colon(\C^n,0)\to(\C^n,0)$ be a bi-Lipschitz homeomorphism such that $\varphi(V(f))=V(g)$. If every irreducible component of the tangent cone $C(V(f),0)$ has a singular locus of dimension $\leq2$, then $m(V(f),0)=m(V(g),0)$.
\end{theorem}
\begin{proof}
Let $\Psi$, $X_i$ and $Y_i$ be as in Theorem \ref{thm:tangent_cones}. Each $X_i$ (resp. $Y_i$) is the zero set of an irreducible homogeneous polynomial, since $C(V(f),0)=V(in(f))$ (resp. $C(V(g),0)=V(in(g))$), and $\Psi$ is a homeomorphism of $\C^n$ with $\Psi(0)=0$ and $\Psi(X_i)=Y_i$. Since $\dim\Sing{X_i}\leq2$, Theorem \ref{thm:top} gives ${\rm deg}(X_i)={\rm deg}(Y_i)$ for all $i$, and, by \eqref{eq:mult_relative} and Theorem \ref{thm:tangent_cones},
$$
m(V(f),0)=\sum_i k_{V(f)}(X_i)\,{\rm deg}(X_i)=\sum_i k_{V(g)}(Y_i)\,{\rm deg}(Y_i)=m(V(g),0).
$$
\end{proof}

\begin{corollary}\label{cor:n4}
Conjecture \ref{conj:metric} has a positive answer for $n\leq4$.
\end{corollary}
\begin{proof}
Let $f,g\colon(\C^n,0)\to(\C,0)$, $n\leq4$, be germs of reduced holomorphic functions, and assume that there is a bi-Lipschitz homeomorphism $\varphi\colon(\C^n,0)\to(\C^n,0)$ such that $\varphi(V(f))=V(g)$. Each irreducible component of $C(V(f),0)$ is a reduced cone of dimension $n-1\leq3$, so its singular locus has dimension $\leq n-2\leq2$, and Theorem \ref{thm:main} applies.
\end{proof}

\begin{remark}\label{rem:final}
(a) The proof of Theorem \ref{thm:main} gives the following more general statement: $m(V(f),0)=m(V(g),0)$ if every irreducible component $X_i=V(P_i)$ of $C(V(f),0)$ satisfies $\dim\Sing{X_i}\leq2$ or $\chi(F_{P_i})\neq0$ (by Theorem \ref{thm:top} and Proposition \ref{prop:lefschetz}(a)).

(b) The same argument gives the following topological statement: let $f$ and $g$ be as in Question \ref{question:zariski}; if there is a homeomorphism $\Psi\colon\C^n\to\C^n$ with $\Psi(0)=0$ and $\Psi(C(V(f),0))=C(V(g),0)$ that maps each irreducible component of $C(V(f),0)$ onto an irreducible component of $C(V(g),0)$ with the same relative multiplicity, and if the irreducible components of $C(V(f),0)$ have singular loci of dimension $\leq2$, then $m(V(f),0)=m(V(g),0)$.
\end{remark}

\section{Examples}\label{section:examples}

\begin{example}\label{exam:umbrella}
Let $P=x^2z-y^2w$ on $\C^4$. Since $P=x^2\cdot z-y^2\cdot w$ is of degree one in $(z,w)$ with coprime coefficients $x^2$ and $y^2$, $P$ is irreducible. The partial derivatives of $P$ are $2xz$, $-2yw$, $x^2$ and $-y^2$, so $\Sing{V(P)}=\{x=y=0\}$, a plane, and $M=\mathbb{P} V(P)$ is a cubic surface in $\C\mathbb{P}^3$ which is singular along the line $L=\{x=y=0\}$. Note that $M$ also contains the line $L'=\{z=w=0\}$, which is disjoint from $L$.

Let us compute $\chi(M)$. The map $\rho\colon M\setminus L'\to\C\mathbb{P}^1$, $(x:y:z:w)\mapsto(z:w)$, is algebraic, and its fibre over $(s:t)$ is $\{(x:y:s:t): sx^2=ty^2\}$, which is isomorphic to the curve $\{(x,y)\in\C^2: sx^2=ty^2\}$, a union of at most two lines through the origin; hence it is contractible. By the Fubini formula, $\chi(M\setminus L')=\chi(\C\mathbb{P}^1)=2$, so, by additivity, $\chi(M)=\chi(L')+\chi(M\setminus L')=4$, and, by the proof of Proposition \ref{prop:lefschetz}(b), $\chi(F_P)=3\,(4-\chi(M))=0$. Thus Proposition \ref{prop:lefschetz}(a) does not apply to $V(P)$, and, since $\dim\Sing{V(P)}=2$, neither does \cite[Theorem 6.6]{FernandesS:2023}; Theorem \ref{thm:top} shows that if $Q$ is an irreducible homogeneous polynomial on $\C^4$ and $\Psi\colon\C^4\to\C^4$ is a homeomorphism with $\Psi(0)=0$ and $\Psi(V(P))=V(Q)$, then ${\rm deg}(Q)=3$.

In the notation of Section \ref{section:cones}, $k=1$, $D_1=L$ and $\delta_1=1$. At a point $(0:0:1:a)$, $a\neq0$, of $L$, in the affine chart $z=1$, $M$ is given by $x^2-y^2w=0$ with $w$ close to $a\neq0$, so its transversal type is $A_1$ and, by Lemma \ref{lemma:transversal}, $c_1=1+(-1)^{N-1}\mu_1=0$ ($N=2$, $\mu_1=1$). At the points $(0:0:1:0)$ and $(0:0:0:1)$, $M$ has Whitney umbrella singularities $x^2=y^2w$ and $y^2=x^2z$ (in the charts $z=1$ and $w=1$, respectively). The Milnor fibre of the weighted homogeneous polynomial $y^2-x^2z$ is homotopy equivalent to the global affine fibre $\{y^2-x^2z=1\}$ \cite[\S 9]{Milnor:1968}, which is a double covering of $\C^2_{(x,z)}$ branched along the curve $\{x^2z=-1\}\cong\C^*$; so, by the Fubini formula, its Euler characteristic is $2\cdot\chi(\C^2\setminus\{x^2z=-1\})+\chi(\C^*)=2\cdot1+0=2$. Hence $\Sigma_0=\{(0:0:1:0),(0:0:0:1)\}$, with $c_M=2$ at both points, and $\chi(D_1^\circ)=0$. Since $\chi_{4,3}=9$, formula \eqref{eq:formula} gives $\chi(M)=9-(-1)(0-3)-2=4$, as it should.
\end{example}

\begin{example}\label{exam:dim5}
Let $P=u(x^2+z^2)+v(y^2+z^2)$ on $\C^5$. As above, $P$ is irreducible, since $x^2+z^2=(x+iz)(x-iz)$ and $y^2+z^2=(y+iz)(y-iz)$ are coprime. The partial derivatives of $P$ are $x^2+z^2$, $y^2+z^2$, $2ux$, $2vy$ and $2z(u+v)$. If they all vanish and $z=0$, then $x=y=0$; if they all vanish and $z\neq0$, then $x,y\neq0$, so $u=v=0$, and $x=\pm iz$, $y=\pm iz$. Hence $\Sing{V(P)}$ is the union of the plane $\{x=y=z=0\}$ and of the four lines $\{u=v=0,\ x=\pm iz,\ y=\pm iz\}$, so $\dim\Sing{V(P)}=2$. The plane $\Pi=\{u=v=0\}\cong\C\mathbb{P}^2$ is contained in $M=\mathbb{P} V(P)\subset\C\mathbb{P}^4$, and the map $(u:v:x:y:z)\mapsto(u:v)$ is defined on $M\setminus\Pi$, with fibres isomorphic to the affine quadric cones $\{s(x^2+z^2)+t(y^2+z^2)=0\}\subset\C^3$, which are contractible. Hence, by the Fubini formula and additivity, $\chi(M)=\chi(\Pi)+\chi(\C\mathbb{P}^1)=3+2=5$, and $\chi(F_P)=3\,(5-5)=0$. So Theorem \ref{thm:main} applies to germs in $\C^5$ whose tangent cone is $V(P)$ (e.g., $f=P+h$ with $h\in\mathfrak m^4$), while Proposition \ref{prop:lefschetz}(a) and \cite[Theorem 6.6]{FernandesS:2023} do not.
\end{example}

\bigskip
\noindent {\bf Acknowledgments:} During the preparation of this work, the author used the AI assistants Claude and Gemini to refine arguments, locate references, and draft and revise the manuscript. All AI-generated suggestions were independently verified by the author, who takes full responsibility for the content and accuracy of the final manuscript.

\end{document}